\documentclass [11pt,a4paper]{article}
\usepackage{amsmath}
\usepackage{amsthm}
\usepackage{amssymb}
\usepackage{amscd}
\usepackage{amsfonts}
\usepackage{amsbsy}
\usepackage{epsfig,afterpage}
\usepackage[dvips]{psfrag}
\usepackage{graphicx}
\usepackage{indentfirst, latexsym, bm,amssymb}
\usepackage{bbding}

\advance\textwidth by +1.0in \advance\textheight by +1.0in
\advance\oddsidemargin by -0.5in \advance\evensidemargin by -1.0in
\advance\topmargin by -0.5in
\bibliography{bibfile}

\newtheorem {theorem} {Theorem}%[section]
\newtheorem {proposition} [theorem]{Proposition}
\newtheorem {corollary} [theorem]{Corollary}

\newtheorem {example} [theorem]{Example}

\def\p{\partial}

\title{\large\bf The $16$th Hilbert problem on algebraic limit cycles}
\author{\normalsize\bf\sc Xiang Zhang\footnote{\small
The author is partially supported by NNSF of China grant 10831003, and Shanghai Pujiang Program grant 09PJD013.}\\
\normalsize\it Department of Mathematics, Shanghai Jiaotong
               University,  \\ \normalsize\it Shanghai 200240,
               The People's Republic of China.
      \\ \normalsize  E-mail: xzhang@sjtu.edu.cn\\}
\date{}

\begin{document}
\maketitle

\begin{abstract}
\noindent For real planar polynomial differential systems there
appeared a simple version of the $16$th Hilbert problem on
algebraic limit cycles: {\it Is there an upper bound on the number
of algebraic limit cycles of all polynomial vector fields of
degree $m$?} In [J. Differential Equations, 248(2010), 1401--1409]
Llibre, Ram\'irez and Sadovskia solved the problem, providing an
exact upper bound, in the case of invariant algebraic curves
generic for the vector fields, and they posed the following
conjecture: {\it Is $1+(m-1)(m-2)/2$ the maximal number of
algebraic limit cycles that a polynomial vector field of degree
$m$ can have?}

In this paper we will prove this conjecture for planar polynomial
vector fields having only nodal invariant algebraic curves. This
result includes the Llibre {\it et al}\,'s as a special one. For
the polynomial vector fields having only non--dicritical invariant
algebraic curves we answer the simple version of the 16th Hilbert
problem.

\hskip0.000011mm

\noindent {\bf Key words and phrases:} polynomial differential
systems, holomorphic singular foliations, simple version of the
16th Hilbert problem, algebraic limit cycles.

\hskip0.000011mm

\end{abstract}

\bigskip

\section{Introduction and the statement of the main results}\label{s1}
\setcounter{section}{1} \setcounter{equation}{0}

The second part of the $16$th Hilbert problem still remain open (see for example, \cite{Hi,Sm}), even through some nice results on the upper bounds of the number of limit cycles can be found in the references (see for instance \cite{BNY, CLi,Il, IS} and the references therein). Related to algebraic limit cycles of real planar polynomial vector fields there appeared a simple version of the $16$th Hilbert problem, see Llibre {\it et al} \cite{LRS}.

\noindent{\bf A simple version of the $16$th Hilbert problem:} {\it Is there
an upper bound on the number of algebraic limit cycles of all real planar
polynomial vector fields of a given degree?}

This simple version of the 16th Hilbert problem provides a nice
connection between the two parts of the 16th Hilbert problem.

Consider real planar polynomial vector fields of degree $m$
\begin{equation}\label{ev1}
\mathcal X=p(x,y)\frac{\partial }{\p x}+q(x,y)\frac{\p }{\p y},
\end{equation}
with $p(x,y),q(x,y)\in\mathbb R[x,y]$ the ring of real polynomials
in $x,y$ and $\mbox{max}\{\deg p,\deg q\}=m$, or the associated
polynomial differential systems
\[
\dot x=p(x,y),\qquad \dot y=q(x,y).
\]
An algebraic curve $f=0$ with $f\in\mathbb C[x,y]$ the ring of
polynomials in $x,y$ with coefficients in $\mathbb C$ is {\it
invariant} by the vector field $\mathcal X$ if there exists some
$K\in \mathbb C[x,y]$ such that
\[
\mathcal X f=p\frac{\p f}{\p x}+q\frac{\p f}{\p y} =Kf.
\]
The polynomial $K$ is called the {\it cofactor} of $f$. An {\it
algebraic limit cycle} is a limit cycle which is contained in an
invariant algebraic curve of $\mathcal X$. A {\it limit cycle} of
an analytic vector field is an isolated periodic orbit in the set
of all periodic orbits of the vector field.

The simple version of the 16th Hilbert problem, i.e. the problem
on the upper bound of the number of algebraic limit cycles, is
solved in \cite{LRS} for all real planar polynomial vector fields
which have only irreducible invariant algebraic curves generic. A
set of invariant algebraic curves, saying $f_j=0$, $j=1,\ldots,k$,
of a planar polynomial vector field is {\it generic} if the
following five conditions hold:
\begin{itemize}

\item{} All the curve $f_j=0$ are {\it non--singular}, (i.e. there
are no points of $f_j=0$ at which $f_j$ and its first derivative
all vanish).

\item{} The highest order homogeneous terms of $f_j$ have no
repeated factors.

\item{} If two curves intersect at a point in the affine plane,
they are transversal at this point.

\item{} There are no more than two curves $f_j=0$ meeting at any
point in the affine plane.

\item{} There are no two curves having a common factor in the
highest order homogeneous terms.
\end{itemize}
The main result of Llibre {\it et al} \cite{LRS} proved that {\it
for a real planar polynomial vector field of degree $m$ having all its
irreducible invariant algebraic curves generic, the maximal number
of algebraic limit cycles is at most $1+(m-1)(m-2)/2$ if $m$ is
even, and $(m-1)(m-2)/2$ if $m$ is odd, and the upper bounds can
be reached.} In the same paper the authors' conjecture 3 stated
that

\noindent{\bf Conjecture.} Is $1+(m-1)(m-2)/2$ the maximal number
of algebraic limit cycles that a polynomial vector field of degree
$m$ can have?

Our first result verifies this conjecture for real planar
polynomial vector fields having only nodal invariant algebraic
curves. We say that an algebraic curve $S$ (not necessary
irreducible) is {\it nodal} if all its singularities are of normal
crossing type, that is at any singularity of $S$ there are exactly
two branches of $S$ which intersect transversally.

\begin{theorem}\label{t0}
If a real planar polynomial vector field \eqref{ev1} of degree $m$
has only nodal invariant algebraic curves taking into account the
line at infinity, then the following hold.
\begin{itemize}

\item[$(a)$] The maximal number of algebraic limit cycles of the
vector fields is at most $1+(m-1)(m-2)/2$ when $m$ is even, and
$(m-1)(m-2)/2$ when $m$ is odd.

\item[$(b)$] There exist systems of form \eqref{ev1} which have
the maximal number of algebraic limit cycles.
\end{itemize}
\end{theorem}

We mention that our result is an essential improvement of that given
in \cite{LRS}, because our assumptions only satisfy the third and
fourth conditions of the generic conditions of Theorem $1$ of \cite{LRS}.

Recently Llibre {\it et al} \cite{LRS1} obtained an upper bound on the number of algebraic
limit cycles for real planar polynomial vector fields which have only non--singular
invariant algebraic curves. The main result states that {\it for a
real planar polynomial vector field of degree $n$ having all its irreducible
invariant algebraic curves non--singular, the maximal number of
algebraic limit cycles is at most $n^4/4+3n^2/4+1$}.

We note that the results given in \cite{LRS,LRS1} both require a
sufficient condition that all the invariant algebraic curves of a
prescribed vector field are non--singular, and so they cannot be
self--intersected.

Our next result will study the case that the invariant algebraic
curves may be singular and the vector field has a more general
form than that given in \eqref{ev1}, i.e.
\begin{equation}\label{e1}
\mathcal X=(p(x,y)+xr(x,y))\frac{\partial }{\partial
x}+(q(x,y)+yr(x,y))\frac{\partial }{\partial y},
\end{equation}
where $p,q,r\in\mathbb R[x,y]$, $\max\{\deg p,\deg q, \deg r\}=m$
and $r$ is a homogeneous polynomial or is identically zero. We
also call $m$ the degree of the vector field \eqref{e1}. In the
next section we will give more explanation on the degree $m$. For
people working in real planar polynomial vector fields they
usually call \eqref{e1} a vector field of degree $m+1$ if
$r(x,y)\not\equiv 0$.

Recall that Theorem \ref{t0} has the restriction on the
singularities of the invariant algebraic curves. We now turn to
the case having some assumption on singularities of the vector
fields. We assume that the singularities of the vector field on
the invariant algebraic curves are non--dicritical. A singularity
of a vector field is {\it non--dicritical} if there are only
finitely many invariant integral curves passing through it. An
invariant algebraic curve is {\it non--dicritical} if there is no
dicritical singularities on it. Clearly a non--dicritical
algebraic curve can be singular.

The following is our second main result.

\begin{theorem}\label{t1}
If a real planar polynomial vector field \eqref{e1} of degree $m$
has all its invariant algebraic curves non-dicritical, then the
following hold.
\begin{itemize}

\item[$(a)$] If $r(x,y)\equiv 0$, the maximal number of algebraic
limit cycles of the vector fields is at most $1+m(m-1)/2$ when $m$
is even, and $m(m-1)/2$ when $m$ is odd.

\item[$(b)$] If $r(x,y)\not\equiv 0$, the maximal number of
algebraic limit cycles of the vector fields is at most
$1+(m+1)m/2$ when $m$ is even, and $(m+1)m/2$ when $m$ is odd.
\end{itemize}
\end{theorem}

We note that Theorem \ref{t1} solves the simple version of the
16th Hilbert problem on algebraic limit cycles for real planar
polynomial vector fields having only non--dicritical invariant
algebraic curves. From the proof of this theorem we guess the
upper bound is not the best one. We conjecture that {\it the best
upper bound for the number of algebraic limit cycles in the
non--dicritical case should be the same as that of Theorem
\ref{t0}}. We remark that the invariant algebraic curves in Theorem \ref{t1} may not satisfy
any one of the conditions that the generic algebraic curves have.

Theorem \ref{t1} has an easy consequence.

\begin{corollary}\label{c1}
If a real planar polynomial vector field \eqref{e1} of degree $m$
has no dicritical singularities, then the following hold.
\begin{itemize}

\item[$(a)$] If $r(x,y)\equiv 0$, the maximal number of algebraic
limit cycles of the vector fields is at most $1+m(m-1)/2$ when $m$
is even, and $m(m-1)/2$ when $m$ is odd.

\item[$(b)$] If $r(x,y)\not\equiv 0$, the maximal number of
algebraic limit cycles of the vector fields is at most
$1+(m+1)m/2$ when $m$ is even, and $(m+1)m/2$ when $m$ is odd.
\end{itemize}
\end{corollary}

The following result provides an exact upper bound on the number
of algebraic limit cycles for polynomial vector fields in
the non--dicritical case with an extra assumption.

\begin{theorem}\label{t11}
For real planar polynomial vector fields \eqref{e1} of degree
$m\ge 2$ having no dicritical singularities, if they have at least
three invariant algebraic curves then the following hold.
\begin{itemize}

\item[$(a)$] The maximal number of algebraic limit cycles of the
vector fields is at most $1+(m-1)(m-2)/2$ when $m$ is even, and
$(m-1)(m-2)/2$ when $m$ is odd.

\item[$(b)$] The maximal number can be reached only for some
polynomial vector fields \eqref{e1} with  $r(x,y)\not\equiv 0$ and
the number of invariant algebraic curves to be three.
\end{itemize}
\end{theorem}

Theorem \ref{t11} has verified Conjecture 3 of \cite{LRS} in the
non--dicritical case with the extra assumption on the number of
invariant algebraic curves. Its proof follows from those of
Theorems \ref{t0} and \ref{t1}, the details are omitted.

This  paper is organized as follows. In the next section we will
present some backgrounds on the degree of invariant algebraic
curves for holomorphic singular foliations. In Section \ref{s3} we
will prove our main results. The last section is an appendix, which provides a
proof to Proposition \ref{p1}.

\section{Upper bound on the degree of invariant algebraic curves
}\label{s2} \setcounter{section}{2}
\setcounter{equation}{0}\setcounter{theorem}{0}

Let $\cal F$ be a holomorphic singular foliation by curves of the
complex projective plane $\mathbb {CP}$(2). Taking an affine
coordinate system $(x,y)$ such that $\mathcal F$ are the solutions
of $\tilde Pdy-\tilde Qdx=0$. Let $L$ be a straight line which is
not invariant by $\mathcal F$. Then the maximal number of the
points $p\in L$ such that either $p\in\{(x,y);\,
P(x,y)=Q(x,y)=0\}$ or the leaf of $\mathcal F$ through $p$ is
tangent to $L$ is bounded by $\max\{\deg P,\deg Q\}$. For a
generical line $L$, this maximal number is a constant. We call it
the degree of $\mathcal F$.

Consider a holomorphic singular foliation $\mathcal F$ of degree
$m$. In the projective coordinates, $\cal F$ can be written as the
closed one--form
\[
\widetilde\omega=P(X,Y,Z)dX+Q(X,Y,Z)dY+R(X,Y,Z)dZ,
\]
where $P,Q,R\in \mathbb C[X,Y,Z]$ are homogeneous polynomials of
degree $m+1$ satisfying the projective condition $XP+YQ+ZR=0$. As
usual, $\mathbb C[X,Y,Z]$ denotes the complex polynomial ring in
the homogeneous coordinates $X,Y$ and $Z$. In the affine
coordinates, $\cal F$ can be written as the one--form
\[
\omega=-(q(x,y)+yr(x,y))dx+(p(x,y)+xr(x,y))dy,
\]
or as the vector field
\[
{\cal X}=(p(x,y)+xr(x,y))\frac{\partial }{\partial
x}+(q(x,y)+yr(x,y))\frac{\partial }{\partial y},
\]
where $p,q,r\in \mathbb C[x,y]$ with $\max\{\deg p,\deg q,\deg
r\}=m$ and $r(x,y)$ is a homogeneous polynomial of degree $m$ or
is naught. If $r\equiv 0$ then $\max\{\deg p,\deg q\}=m$. These
claims can be found in \cite {Li} and \cite {CL}.

A point $(X_0,Y_0,Z_0)\in \mathbb CP(2)$ is called a {\it
singularity} of $\cal F$ if $P(X_0,Y_0,Z_0) =Q(X_0,Y_0,Z_0)$ $
=R(X_0,Y_0,Z_0)=0$; or in affine plane $(X_0,Y_0,Z_0)=(x_0,y_0,1)$
satisfies $p(x_0,y_0)+x_0r(x_0,y_0)=q(x_0,y_0)+y_0r(x_0,y_0)=0$. A
singularity of $\cal F$ is called {\it non--dicritical} if there
are only finitely many integral curves passing through it.
Otherwise, it is called {\it dicritical}.

An algebraic curve $S$ defined by a reduced homogeneous polynomial
$F(X,Y,Z)\in \mathbb C[X,Y,Z]$ is called {\it invariant} by $\cal
F$ if $\widetilde \omega\wedge dF=F\theta$, where $\theta$ is a
two--form. Recall that a {\it reduced polynomial} is the one which
has no repeat factors. In what follows, for simplicity we also say
$F$ is an invariant algebraic curve. It is easy to prove \cite{Zh}
that $F$ is an invariant algebraic curve if and only if ${\cal
X}f=kf$ for some $k\in \mathbb C[x,y]$, where
$f=\left.F\right|_{Z=1}$.

Theorem 1 of Cerveau and Lins Neto \cite{CL} in 1991 obtained the
exact upper bound on the degree of nodal invariant algebraic
curves, which is the key point to prove Theorem \ref{t0}.

\begin{theorem}\label{CLt} {\rm (Cerveau and Lins Neto 1991)}
Let $\mathcal F$ be a foliation in $\mathbb CP(2)$ of degree $m$,
having $S$ as a nodal invariant algebraic curve with the reduced
homogeneous equation $F=0$ of degree $n$. Then $n\le m+2$.
Moreover if $n=m+2$ then $F$ is reducible and the foliation
$\mathcal F$ is of logarithmic type, that is given by a rational
closed form $\sum\limits_{i}\lambda_i\frac{dF_i}{F_i}$, where
$\lambda_i\in \mathbb C$ and $F_i$ are the irreducible homogeneous
components of $F$ and $\sum\limits_{i}\lambda_i\deg F_i=0$.
\end{theorem}

In the non--dicritical case  Carnicer \cite{Ca} in 1994 obtained
the same upper bound as that given in Theorem \ref{CLt}, which
solved the Poincar\'{e} problem \cite{Po} in the non--dicritical
case. We will use it to prove our Theorem \ref{t1}.

\begin{theorem} \label{t2} {\rm (Carnicer 1994)}  Let $\cal F$ be a holomorphic
singular foliation of degree $m$ in $\mathbb CP(2)$. Assume that
$S$ is an algebraic curve which is invariant by $\cal F$,  and is
given by a reduced polynomial $F$ of degree $n$. If there are no
dicritical singularities of $\cal F$ on $S$, then $n\le m+2$.
\end{theorem}

In the proof of this last result, the author had used the following result, which is due to Cerveau and Lins Neto \cite{CL}.

\begin{proposition} \label{p1}  Assume that $\cal F$ is a
holomorphic singular foliation of degree m in $\mathbb {CP}(2)$,
and that $S$ is a reduced algebraic curve of degree $n$ which is
invariant by $\cal F$. Let $\chi(S)$ be the {\it intrinsic Euler characteristic} of $S$
$($see {\rm\cite{Li}}$)$ and let $g(S)$ be the topological genus of $S$.
Then
\begin{equation}\label{sb2}
\chi(S)=2-2g(S)=\sum_{B}\mu_p({\cal F},B)-n(m-1),
\end{equation}
where the sum runs over all the local branches $B$ of $S$ passing
through the singularities of $\cal F$ in $S$, and $\mu_p({\cal F},B)$ is the multiplicity of $\cal F$ at $B$ passing through the singularity $p$.
\end{proposition}

Since the proof of the last result given in \cite{CL} has a gap inside, we will present a new proof to it in the appendix. The multiplicity of $\cal F$ at $B$ passing through $p$ is defined as follows:
for each singularity $p$ of $\cal F$ such that $p\in S$, and each
local branch $B$ of $S$ passing through $p$, take a vector field
${\cal X}=P\frac{\partial}{\partial x}+Q\frac{\partial}{\partial
y}$ to represent $\cal F$ in a neighborhood of $p$ and a minimal
Puiseux's parameterization of $B$, saying that
$\phi=(\phi_1,\phi_2):\,\mathbb D\rightarrow\mathbb C^2$ such that
$\phi(0)=0$, where $\mathbb D$ is a disk centered at $0\in
\mathbb C$. We define the {\it multiplicity} of $\cal F$ at $B$
to be the order of $\phi^*({\cal X})=R(t)\frac{d}{dt}$ at $t=0\in
\mathbb D$, denoted by $\mu_p({\cal F},B)$. Then
\[
\mu_p({\cal F},B)=\frac{1}{2\pi
i}\int_{\gamma(B)}\frac{dR(t)}{R(t)},
\]
where $\gamma(B)=re^{i\theta}$, $r>0$ small, is the homology class
in $H_1(B\setminus\{p\},t)$ of the curve
$\theta\rightarrow\phi(re^{i\theta})$, $0\le\theta\le 2\pi$.

As a by--product of the equality \eqref{sb2} we have the following well--known result. Since the proof is short, we will present it in the appendix.

\begin{corollary}\label{c2}
An irreducible non--singular algebraic curve $S$ of degree $n$ has
the Euler characteristic $\chi(S)=-n(n-3)$.
\end{corollary}

\section{Proof of the main results}\label{s3}
\setcounter{section}{3}
\setcounter{equation}{0}\setcounter{theorem}{0}

For proving the theorems we need the following Harnack's theorem,
for a proof see for instance \cite{Co,Vi,Wi}.

\begin{theorem}\label{t3} {\rm (Harnack's Theorem)}
The number of ovals of a real irreducible algebraic curve of
degree $n$ is at most
\[
1+(n-1)(n-2)/2-\sum\limits_p\nu_p(S)(\nu_p(S)-1),
\]
if $n$ is even, or
\[
(n-1)(n-2)/2-\sum\limits_p\nu_p(S)(\nu_p(S)-1),
\]
if $n$ is odd, where $p$ runs over all the singularities of
$\mathcal F$ on $S$, and $\nu_p(S)$ is the order of $S$ at the singular
point $p$. Moreover these upper bounds can be reached
for convenient algebraic curves of degree $n$.
\end{theorem}

The following result, due to Giacomini, Llibre and Viano
\cite{GLV}, provides the location of limit cycles for a real
planar differential system having an inverse integrating factor,
for a different proof see \cite{LR}.

\begin{theorem}\label{tGLV} Let $\mathcal X$ be a $C^1$ vector
field defined in the open subset $U$ of $\mathbb R^2$, and let $V: \,
U\rightarrow \mathbb R$ be an inverse integrating factor of
$\mathcal X$. If $\gamma$ is a limit cycle of $\mathcal X$, then
$\gamma$ is contained in $\{(x,y)\in U:\, V(x,y)=0\}$.
\end{theorem}

\subsection{Proof of Theorem \ref{t0}}
\noindent $(a)$.  Write system \eqref{ev1} in the one--form
\[
q(x,y)dx-p(x,y)dy.
\]
Its projective one--form is
\begin{equation}\label{e01}
\omega_0=ZQdX-ZPdY+(YP-XQ)dZ,
\end{equation}
where $X,Y,Z$ are the homogeneous coordinates and
\[
P=Z^mp(X/Z,Y/Z),\,\,\,\,\, Q=Z^mq(X/Z,Y/Z).
\]

Consider the holomorphic singular foliation $\mathcal F_0$ induced
by the one--form $\omega_0$. Clearly $\mathcal F_0$ has the
infinity as an invariant line. Under the assumption of the
theorem, we get from Theorem \ref{CLt} that the total degree $n$
of all invariant algebraic curves of the foliation $\mathcal F_0$
is no more than $m+2$.

\noindent {\it Case} 1.  $n=m+2$. Theorem \ref{CLt} shows that $F$
is reducible, saying $F=F_1\cdot\ldots\cdot F_k$ the irreducible
decomposition with $k\ge 2$. The one--form $\omega_0$ has the
expresssion
\[
\omega_0=F\sum\limits_{i=1}\limits^{k}\lambda_i\frac{dF_i}{F_i},
\]
where $\lambda_i\in \mathbb C$. Recall that an invariant algebraic
curve of a real system can be complex. If it happens its conjugate
is also an invariant algebraic curve of the system. The one--form
$\omega_0$ has the inverse integrating factor $F$, and
consequently is Darboux integrable with the Darboux first integral
$H(X,Y,Z)=F_1^{\lambda_1}\cdot\ldots\cdot F_k^{\lambda_k}$. For
more information on the Darboux theory of integrability, see for
instance \cite{Ll,LZ,LZ1}.

Since the one--form $\omega_0$ is projective, i.e.
$i_E\omega_0=0$, where $E=X\frac{\p}{\p X}+Y\frac{\p}{\p
Y}+Z\frac{\p}{\p Z}$ is the radial vector field and $i_E$ is the
interior productor by $E$, we should have
$\lambda_1\deg{F_1}+\ldots+\lambda_k\deg{F_k}=0$.

If $k=2$, the
foliation $\mathcal F_0$ has a rational first integral
\[
H(X,Y,Z)=F_1^kF_2^{-l},\qquad \mbox{ with } k,l\in\mathbb N,\,\,\,
(k,l)=1, \mbox{ and } k/l=\deg F_2/\deg F_1.
\]
In this case there are infinitely many invariant algebraic curves.
Of course they are not possible of nodal type. Otherwise it is in
contradiction with Theorem \ref{CLt}. So we must have $k\ge 3$.

For $k\ge 3$, we get from the Harnack's theorem that each
invariant algebraic curve has at most $(\deg F_i-1)(\deg
F_i-2)/2+a_i$ ovals, where $a_i=1$ if $\deg F_i$ is even, and
$a_i=0$ if $\deg F_i$ is odd. So the total number of ovals
contained in $F_i$ for $i=1,\ldots,k$ is no more than
\[
\sum\limits_{i=1}\limits^k\left(\frac{(\deg F_i-1)(\deg
F_i-2)}{2}+a_i\right)\le
\frac{(m+2-k)(m+1-k)}{2}+\sum\limits_{i=1}\limits^ka_j,
\]
where we have used Lemma 6 of \cite{LRS} and $\deg F_1+\ldots+\deg
F_k=m+2$. Furthermore the equality holds if and only if one of the
$F_i$'$s$ has the degree $m+3-k$ and the others all have degree
$1$.

Set
\[
M(k)=\frac{(m+2-k)(m+1-k)}{2}+\sum\limits_{j=1}\limits^ka_j.
\]
Then the maximum of the $M(k)$ for $k\in\{3,\ldots,m+2\}$ takes
place when $k=3$, because $\sum\limits_{j=1}\limits^ka_j\le
[m/2]+1$, where $[\cdot]$ denotes the integer part function. For
$k=3$ and the three invariant algebraic curves have respectively
the degrees 1,1 and $m$, the maximum is
\[
\frac{(m-1)(m-2)}{2}+a,
\]
where $a=1$ if $m$ is even and $a=0$ if $m$ is odd.

\noindent {\it Case} 2. $n\le m+1$. Recall that the line at
infinity is invariant by the foliation $\mathcal F_0$. If $n\le m$
then the total degree of the invariant algebraic curves in the
affine plane is less than $m$. By the Harnack's theorem we get
from the proof of case 1 that the number of algebraic limit cycles
is less than the maximal value.

If $n=m+1$, the total degree of the invariant algebraic curves in
the affine plane is $m$. By the Harnack's theorem the number of
algebraic limit cycles is less than or equal to the maximal value.
This proves statement $(a)$.

\noindent $(b)$. We only need to prove that  there exists a real
planar polynomial system of form \eqref{ev1} with degree $m$ which
has the maximal number of algebraic limit cycles and the total
degree of the invariant algebraic curves in the affine plane is
$m$ and $m+1$ respectively, because the line at infinity is
invariant.

\noindent{\it Case} 1. The number $m+1$ is the total degree of the invariant
algebraic curves in the affine plane. By the Harnack's theorem
there exists a nonsingular algebraic curve of degree $m$ which has
the maximal number, i.e. $(m-1)(m-2)/2+a$, of ovals, where $a=1$
if $m$ is even, or $a=0$ if $m$ is odd. Denote by $f_1$ this
curve. Choose a straight line, called $f_2$, as the line at
infinity in such a way that which is outside the ovals of $f_1$
and intersects $f_1$ transversally. Choose another straight line,
called $f_3$, which is outside the ovals of $f_1$ and intersects
$f_1$ and $f_2$ transversally and does not meet the intersection
points of $f_1$ and $f_2$.

Let $F_1$, $F_2$ and $F_3$ be the projectivization of $f_1$, $f_2$
and $f_3$, respectively. Taking
$\lambda_1,\lambda_2,\lambda_3\in\mathbb R$ non--zero such that
$\lambda_1m+\lambda_2+\lambda_3=0$ and
$\lambda_i/\lambda_j\not\in\{r\in\mathbb Q;\, r<0\}$. Then the
foliation $\mathcal F_m$ induced by the projective one--form
$\lambda_1F_2F_3dF_1+\lambda_2F_1F_3dF_2+\lambda_3F_1F_2dF_3$ has
only the three invariant algebraic curves $F_1,F_2,F_3$. Hence
$\mathcal F_m$ has exactly $(m-1)(m-2)/2+a$ algebraic limit
cycles. In fact $\mathcal F_m$ has the inverse integrating factor
$F_1F_2F_3$. By Theorem \ref{tGLV} $\mathcal F_m$ has no other
limit cycles, i.e. the non--algebraic ones.

We note that $\mathcal F_m$ is a holomorphic singular foliation of
degree $m$. Since it has the line at infinity invariant, its
affine expression should be a polynomial differential system of
degree $m$ having the form \eqref{ev1}.

\noindent{\it Case} 2. The number $m$ is the total degree of the invariant
algebraic curves in the affine plane. In fact the proof can be
obtained from \cite{Ch,LRS}. For completeness and because it is
short, we present it here for readers' convenience.

By the Harnack's theorem there exists a nonsingular algebraic
curve of degree $m$ which has the maximal number, i.e.
$(m-1)(m-2)/2+a$, of ovals, where either $a=1$ or $a=0$ if $m$ is
either even or odd. Denote by $g(x,y)$ this nonsingular algebraic
curve. Choose a linear function $h(x,y)$ such that $h=0$ does not
intersect the ovals of $g=0$, and choose $a,b\in\mathbb R$
satisfying $ah_x+bh_y\ne 0$, then the real planar differential
system
\begin{equation}\label{eee}
\dot x=a g-hg_y,\qquad \dot y=bg+hg_x,
\end{equation}
is of degree $m$ and has all the ovals of $g=0$ as hyperbolic
limit cycles. Moreover system \eqref{eee} has no other limit
cycles. This proves statement $(b)$ and consequently the theorem.
$\Box$

\subsection {Proof of Theorem \ref{t1}} Write system \eqref{e1}
in the one--form
\[
(q(x,y)+yr(x,y))dx-(p(x,y)+xr(x,y))dy.
\]
Its projective one--form is
\begin{equation}\label{e02}
\omega_1=(ZQ+YR)dX-(ZP+XR)dY+(YP-XQ)dZ,
\end{equation}
where $X,Y,Z$ are the homogeneous coordinates and
\[
P=Z^mp(X/Z,Y/Z),\,\,\, Q=Z^mq(X/Z,Y/Z),\,\,\, R=Z^mr(X/Z,Y/Z).
\]

Let $\mathcal F_1$ be the holomorphic singular foliation induced
by $\omega_1$. By the assumption of the theorem $\mathcal F_1$ has
all the invariant algebraic curves non--dicritical, and their
total degree is less than or equal to $m+2$ by Theorem \ref{t2}.

\noindent $(a)$ If $r(x,y)\equiv 0$, the line at infinity is
invariant by the foliation $\mathcal F_1$. So it follows from the
proof of Theorem \ref{t0} that the total degree $n$ of all the
invariant algebraic curves in the affine plane is at most $m+1$.
Recall that $m$ is the degree of the polynomial vector field. From
the proof of case 1 of statement $(b)$ of Theorem \ref{t0}, we
know that there is a foliation of degree $m$ which has invariant
algebraic curves with the total degree $m+2$ taking into account
the line at infinity. Of course, it is reducible that the
invariant algebraic curves by the foliation constructed in case 1
of the proof of statement $(b)$ of Theorem \ref{t0}.

If $\mathcal F_1$ has an irreducible invariant algebraic curve of
degree $m+1$ in the affine plane with the maximal number of ovals
that an algebraic curve of degree $m+1$ can have by Theorem
\ref{t3}, then the foliation has the maximal number of algebraic
limit cycles. In all the other cases there is not a system of form
\eqref{e1} which has the maximal number of algebraic limit cycles.
This proves statement $(a)$.

\noindent $(b)$ If $r(x,y)\not\equiv 0$, the line at infinity is
not invariant by $\mathcal F_1$.  We get from Theorem \ref{t2}
that the total degree $n$ of all invariant algebraic curves of
\eqref{e1} in the affine plane is at most $m+2$. We claim that
there exists a system of form \eqref{e1} having degree $m$ with
$r(x,y)\not\equiv 0$ which has invariant algebraic curves of total
degree $m+2$.

We now prove the claim. Let $f_1,\ldots,f_k\in\mathbb C[x,y]$ with
$k\ge 3$ be reduced such that $\deg f_1+\ldots +\deg f_k=m+2$ and
their projective curves in $\mathbb CP(2)$ defined by
$F_1,\ldots,F_k$ the projectivization of $f_1,\ldots,f_k$ are
nonsingular and intersect transversally and no more than two
curves meeting at the same point. Taking
$\lambda_1,\ldots,\lambda_k\in\mathbb C$ non--zero such that
$\lambda_1\deg F_1+\ldots+\lambda_k\deg F_k=0$ and
$\lambda_i/\lambda_j\not\in\{r\in\mathbb Q;\, r<0\}$ for $1\le
i\ne j\le k$. Then the foliation $\mathcal F^*$ induced by the
projective one--form
$\omega^*=\sum\limits_{j=1}\limits^k\lambda_j\left(\prod\limits_{i=1,i\ne
j}\limits^{k} F_i dF_j\right)$ has degree $m$ and has only the
invariant algebraic curves defined by $F_1,\ldots,F_k$.
Furthermore all the singularities of $\mathcal F^*$ are
non--dicritical \cite{Se}, because they are the intersection
points of $F_1,\ldots,F_k$ and the invariant curves passing
through these singularities are only the branches of $F_i$ for
$i=1,\ldots,k$. By Theorem \ref{t2} the total degree of all
invariant algebraic curves by $\mathcal F^*$ is at most $m+2$.
While $f_1,\ldots,f_k$ have the total degree $m+2$. This implies
that the line at infinity of $\mathcal F^*$ is not invariant. So
its affine expression of $\mathcal F^*$ must have the form
\eqref{e1} with $r(x,y)\not\equiv 0$. This proves the claim.

If $\mathcal F_1$ has an irreducible invariant algebraic curve of
degree $m+2$ with the maximal number of ovals that an algebraic
curve of degree $m+2$ can have by Theorem \ref{t3}, then the
foliation has the maximal number of algebraic limit cycles. In all
the other cases there is not a system of form \eqref{e1} which has
the maximal number of algebraic limit cycles.

We complete the proof of the theorem. $\Box$

We mention that the foliation $\mathcal F^*$ of degree $m$
constructed in the proof of statement $(b)$ of Theorem \ref{t1}
has at least three invariant algebraic curves with the total
degree $m+2$. We do not know if there is a holomorphic singular
foliation of degree $m$ which has a non--dicritical irreducible
invariant algebraic curve of degree either $m+1$ or $m+2$. Of
course as shown in Theorem \ref{CLt} it is not possible for nodal
invariant algebraic curves. We guess it is also not possible for
non--dicritical invariant algebraic curves, but we cannot prove it
now.

Finally we provide an easy example showing the foliation $\mathcal
F^*$ mentioned above.

\noindent {\it Example.} For an algebraic curve $S$ in $\mathbb
CP(2)$ which has the affine representation $f=xy(y-x-1)$. The
projective homogeneous form of $f$ is $F=XY(Y-X-Z)$. Then the
holomorphic foliation $\mathcal F_3^*$ given by the one--fom
\[
Y\left(\lambda_1Y+\lambda_2X-\lambda_1Z\right)dX-
X\left(\lambda_1Y+\lambda_2X+\lambda_2Z\right)dY-\lambda_3XYdZ,
\]
has degree 1 and has only the invariant
algebraic curves $F$ provided that $\lambda_1+\lambda_2+\lambda_3=0$, $\lambda_i\ne 0$ for
$i=1,2,3$ and $\lambda_i/\lambda_j$ for $1\le i\ne j\le 3$ non--negative rational numbers.
The line at infinity, i.e. $Z=0$, is not
invariant for $\mathcal F_3^*$. The singularities of $\mathcal F_3^*$ is
non--dicritical, see \cite{Se}.

\section{Appendix}\label{sa}
\setcounter{section}{4}
\setcounter{equation}{0}\setcounter{theorem}{0}

\subsection{Proof of Proposition \ref{p1}}
Since $S$ is an invariant algebraic curve of degree $n$, we can
choose an affine coordinate system $(x,y)$ such that $S$ cuts the
line at infinity $l_\infty$ transversely at exactly $n$ points.
Let ${\cal X}=P(x,y)\frac{\partial}{\partial
x}+Q(x,y)\frac{\partial}{\partial y}$ represent $\cal F$ in this
coordinate system. Without loss of generality, we suppose that
$p=(1:0:0)$ belongs to $S\cap l_\infty$. Making the change of the
variables $u=\frac {y}{x}$, $v=\frac{1}{x}$, the vector field
$\cal X$ becomes
\[
\widetilde {\cal X}=v^{-m+1}\left[\left(-u\widetilde
P(u,v)+\widetilde Q(u,v) \right)\frac{\partial}{\partial
u}-v\widetilde P(u,v)\frac{\partial}{\partial v}\right]
\]
where $\widetilde
P(u,v)=v^mP\left(\frac{1}{v},\frac{u}{v}\right)$, $\widetilde
Q(u,v)=v^mQ\left(\frac{1}{v},\frac{u}{v}\right)$.

In the coordinate system $(u,v)$, since $S$ intersects $l_\infty$
transversely we can take $u=\psi(v)$ as the local branch
$B_\infty$ of $S$ passing through the singularity $(1:0:0)$.
Clearly, $\psi$ is analytic in $v$. Using the change of the
variables $\alpha=u-\psi(v)$, $\beta=v$, the vector field
$\widetilde{\cal X}$ can be written as
\[
{\widetilde {\cal
X}}^*=\beta^{-m+1}\left[\left(-(\alpha+\psi(\beta)){\widetilde
P}^*+{\widetilde Q}^*+\beta\psi^{\prime}(\beta){\widetilde P}^*
\right)\frac{\partial}{\partial \alpha}-\beta{\widetilde
P}^*\frac{\partial}{\partial \beta}\right],
\]
where ${\widetilde P}^*=\widetilde P(\alpha+\psi(\beta),\beta)$,
${\widetilde Q}^*=\widetilde Q(\alpha+\psi(\beta),\beta)$. Since
$h(u,v)=u-\psi(v)$ is an analytic solution of $\widetilde{\cal
X}$, there exists a locally analytic function $k(u,v)$ such that
\[
\left(-u\widetilde P(u,v)+\widetilde Q(u,v) \right)\frac{\partial
h }{\partial u}-v\widetilde P(u,v)\frac{\partial h}{\partial
v}=hk.
\]
Hence, we have
\[
-(\alpha+\psi(\beta)){\widetilde P}^*+{\widetilde
Q}^*+\beta{\widetilde P}^*\psi^{\prime}(\beta)={\widetilde
k}^*\alpha,
\]
where ${\widetilde k}^*=k(\alpha+\psi(\beta),\beta)$. This shows
that on $B_{\infty}$
\[
{\widetilde {\cal X}}^*=\beta^{-m+1}\left({\widetilde k}^*
\alpha\frac{\partial}{\partial \alpha}-\beta{\widetilde
P}^*\frac{\partial}{\partial \beta}\right).
\]
Set ${\widetilde P}^*|_{\alpha=0}=\beta^l\widehat{P}^*(\beta)$
such that $\widehat{P}^*(0)\ne 0$, and set
\[
\zeta=\frac{{\widehat P}^*}{\beta^{m-2-l}}=\frac{{\widehat
P}^*}{|\beta|^{2(m-2-l)}}\overline\beta^{m-2-l},
\]
where $\overline\beta$ denotes the conjugacy of $\beta$. Then the
vector field $\left.{\widetilde {\cal X}}^*\right|_{\alpha=0}$ at
$\beta=0$ has the multiplicity or a pole of order
\[
\frac{1}{2\pi
i}\int_{\gamma}\frac{d\zeta}{\zeta}=\frac{m-2-l}{2\pi
i}\int_{\gamma}\frac{d\overline\beta}{\overline\beta}=-(m-2-l),
\]
where $\gamma$ is the homology class in $H_1(B_\infty,\beta)$ of
the curve $\theta\rightarrow\phi(re^{i\theta})$ on $\alpha=0$,
$0\le\theta\le 2\pi$. Moreover, from the expression of $\widetilde
{\mathcal X}^* $ we can get easily that $\mu_p({\cal
X},B_\infty)=l+1$.

Let $\pi:\, \widetilde S\rightarrow S$ be a resolution of $S$ by
blowing-ups at the singularities of $S$. Then $\widetilde S$ is
smooth and $2-2g(S)=\chi(\widetilde S)$, which is the Euler
characteristic of $\widetilde S$. We define the intrinsic Euler
characteristic  $\chi(S)$ to be $\chi(\widetilde S)$, and the
vector field in $\widetilde S$ associated with $\cal X$ to be $\pi
^*(\left.{\cal X}\right|_S)=\widetilde {\cal X}$. For each
singularity $p$ of $\cal F$ in $S$, and each local branch $B$ of
$S$ passing through $p$, we obtain a singularity $\widetilde p$ of
$\widetilde {\cal X}$ in $\widetilde S$ and a unique local branch
$\widetilde B$ passing through $\widetilde p$ which is invariant
by $\widetilde{\cal X}$. Then the Poincar\'{e}-Hopf's index of
$\widetilde {\cal X}$ with respect to $\widetilde B$ at
$\widetilde p$ is $\mu_p({\cal F},B)$.

From the choice of the local coordinate system at the beginning of
the proof of this proposition, we know that $l_\infty\cap S$
contains $n$ points, denoted by $p_i$, $i=1,\ldots,n$. We denote
by $l_i$ associated to $p_i$ the quantity $l$ in the above proof
for the singularity $p$. Then we get from the Poincar\'{e}-Hopf's
Index Theorem that
\[
\chi(S)=\sum\limits_{B}\mu_p(\mathcal F,
B)-\sum\limits_{i=1}\limits^n(m-2-l_i),
\]
where $B$ is taken over all the local branches of $S$ passing
through the singularities at the finite plane. Since
$\mu_{p_i}({\cal X},B_\infty)=l_i+1$ we have
\[
\chi(S)=\sum\limits_{B}\mu_p(\mathcal F, B)-n(m-1),
\]
where $B$ is taken over all the local branches of $S$ passing
through the singularities.
We complete the proof of the proposition. \hskip 2.8in $\Box$

Next we provide some examples showing the application of
Proposition \ref{p1}.

\begin{example} \label{ex1} Consider the foliation $\mathcal F_1$ of $\mathbb
CP(2)$ given by the homogeneous differential form
\[
\alpha YZdX+\beta XZdY-(\alpha+\beta)XYdZ,
\]
with $\alpha,\beta\in\mathbb C\setminus\{0\}$ and
$\frac{\alpha}{\beta}\not\in\mathbb R$ $($this assures that all
the singularities of $\mathcal F_1$ are non--dicritical$)$. The
line $X=0$ is invariant by the foliation $\mathcal F_1$, on which
there are two singularities: $P_1=(0:1:0)$ and $P_2=(0:0:1)$. The
vector field associated with $\mathcal F_1|_{X=0}$ at $P_1$ is
$\alpha z\frac{\partial}{\partial z}$, so $\mu_{P_1}(\mathcal F_1,
X=0)=1$. Similarly, The vector field associated with $\mathcal
F_1|_{X=0}$ at $P_2$ is $\alpha y\frac{\partial}{\partial y}$, so
$\mu_{P_2}(\mathcal F_1, X=0)=1$. In addition we have
$\chi(X=0)=2$. Since the foliation is of degree $1$, this verifies
the proposition.
\end{example}

\begin{example} \label{ex2} Consider the foliation $\mathcal F_2$ of $\mathbb
CP(2)$ given by the homogeneous differential form
\[
(2YZ-X^2)ZdX+X(Y+Z)ZdY+(X^3-XY^2-3XYZ)dZ.
\]
The foliation $\mathcal F_2$ has $X=0$ as an invariant line, which
contains exactly two singularities of $\mathcal F_2$:
$P_1=(0:1:0)$ and $P_2=(0:0:1)$. We can check easily that $P_1$
and $P_2$ are both non--dicritical. The vector field associated
with $\mathcal F_2|_{X=0}$ at $P_1$ and $P_2$ are $-2
z^2\frac{\partial}{\partial z}$ and $-2 y\frac{\partial}{\partial
y}$, respectively. So we have $\mu_{P_1}(\mathcal F_2, X=0)=2$ and
$\mu_{P_2}(\mathcal F_2, X=0)=1$. Now the foliation has degree
$2$, this verifies the proposition.
\end{example}

We note that $P_1$ and $P_2$ are respectively dicritical and
nondicritical singularities of $\mathcal F_2$.

\begin{example} \label{ex3} Consider the foliation $\mathcal F_3$ of $\mathbb
CP(2)$ given by the homogeneous differential form
\[
(X^3-2Y^2Z)ZdX-X(Y^2+Z^2)ZdY-(X^4-2XY^2Z-XYZ^2-XY^3)dZ.
\]
The foliation $\mathcal F_3$ has also $X=0$ as an invariant line,
on which there are only the non--dicritical singularities:
$P_1=(0:1:0)$ and $P_2=(0:0:1)$.  The vector field associated with
$\mathcal F_3|_{X=0}$ at $P_1$ and $P_2$ are $-2
z^2\frac{\partial}{\partial z}$ and $-2
y^2\frac{\partial}{\partial y}$, respectively. So we have
$\mu_{P_1}(\mathcal F_3, X=0)=2$ and  $\mu_{P_2}(\mathcal F_3,
X=0)=2$. Now the foliation has degree $3$, this verifies the
proposition.
\end{example}

We can check that that $P_1$ and $P_2$ are both dicritical
singularities of $\mathcal F_3$, in fact they are saddle node.

\subsection{Proof of Corollary \ref{c2}}
Taking an affine coordinate system $(x,y)$ of
$\mathbb CP(2)$ such that $S$ intersects the line at infinity
transversally. Denote
by $p_i$ and $B_i$, $i=1,\ldots,n$, the $n$ intersection points
and the $n$ branches of $S$ passing through $p_i$, respectively. Let
$f\in\mathbb C[x,y]$ be a reduced equation of the affine part of
$S$.  We denote by $\mathcal G_f$ the holomorphic foliation by curves of $\mathbb CP(2)$ which
extends the foliation of $\mathbb C^2$ given by $df$. Then $\mathcal G_f$ has degree $n-1$. Applying the formula
\eqref{sb2} to the foliation $\mathcal G_f$, we have
\[
\sum\limits_{i=1}\limits^n\mu_{p_i}(\mathcal G_f,B_i)=\chi(S)+n(n-2).
\]
Since $S$ is nonsingular, it follows that $p_i$ for $i=1,\ldots,n$ are the only
singularities of $\mathcal G_f$, which are located an the intersection of $S$ with $l_\infty$. Moreover we have
$\mu_{p_i}(\mathcal G_f,B_i)=1$ for $i=1,\ldots,n$. This shows that $n=\chi(S)+n(n-2)$, and consequently the corollary follows. $\Box$


\begin{thebibliography}{99}

\bibitem{BNY} G. Binyamini, D. Novikov and S. Yakovenko, On the number of zeros of Abelian integrals:
a constructive solution of the infinitesimal Hilbert sixteenth problem, {\it Invent. Math.} {\bf 181} (2010), 227--289.

\bibitem{Ca} {M. M. Carnicer}, The Poincar\'{e} problem in the
non--dicritical case, {\it Ann. Math.} {\bf 140} (1994),
289--294.

\bibitem{CL} {D. Cerveau and A. Lins Neto}, Holomorphic foliations in
$\mathbb{CP}(2)$ having an invariant algebraic curve, {\it Ann.
Inst. Fourier, Grenoble} {\bf 41} (1991),  883--903.

\bibitem{Ch} {C. Christopher}, Polynomial vector fields with prescribed algebraic limit cycles, {\it Geom.
Dedicata} {\bf 88} (2001), 255--258.

\bibitem{CLi} C. Christopher and Chengzhi Li, {\it Limit cycles of differential equations}, Advanced
Courses in Mathematics, CRM Barcelona, Birkh\"{a}user Verlag, Basel, 2007.

\bibitem{Co} {J. L. Coolidge}, {\it A Treatise on Algebraic Plane Curves}, Dover, New York,
1959.

\bibitem{GLV} {H. Giacomini, J. Llibre and M. Viano}, On the nonexistence, existence and uniqueness of limit cycles,
{\it Nonlinearity} {\bf 9} (1996), 501--516.

\bibitem{Hi} D. Hilbert, Mathematical problems, Reprinted from {\it Bull. Amer. Math. Soc.} {\bf 8} (1902), 437-479,
in {\it Bull. Amer. Math. Soc.} (N.S.) {\bf 37} (2000), 407-436.

\bibitem{Il} Yu. S. Il'yashenko, {\it Finiteness Theorems for Limit Cycles}, Transl. Math. Monogr. {\bf 94}, Amer. Math. Soc., Providence, 1991.

\bibitem{IS} I. Itenberg and E. Shustin, Singular points and limit cycles of planar polynomial vector fields,
{\it Duke Math. J.} {\bf 102} (2000), 1--37.


\bibitem{Li} {A. Lins Neto}, Algebraic solutions of polynomial differential
equations and foliations in dimension two, {\it Lect. Notes
Math.} {\bf 1345} (1988), 192--232.

\bibitem {Ll}
{J. Llibre}, {\it Integrability of polynomial differential
systems}, in Handbook of differential equations, Elsevier,
Amsterdam, 2004, pp.437--532.

\bibitem{LRS} {J. Llibre, R. Ram\'irez and N. Sadovskaia}, On the $16$th Hilbert problem
for algebraic limit cycles, {\it J. Differential Equations} {\bf
248} (2010), 1401--1409.

\bibitem{LRS1} {J. Llibre, R. Ram\'irez and N. Sadovskaia}, On the $16$th Hilbert problem
for limit cycles on nonsingular algebraic curves, {\it J.
Differential Equations} (2011), {\bf 250},  983--999.

\bibitem{LR} {J. Llibre and G. Rodr\'iguez}, Configurations of limit cycles and planar polynomial vector fields,
{\it J. Differential Equations} {\bf 198} (2004), 374--380.

\bibitem{LZ} J. Llibre and Xiang Zhang, Darboux theory of integrability in $\mathbb C^n$
taking into account the multiplicity, {\it J. Differential Equations} {\bf 246} (2009), 541--551.

\bibitem{LZ1} J. Llibre and Xiang Zhang, Rational first integrals in the Darboux theory of
integrability in $\mathbb C^n$, {\it Bull. Sci. Math.} {\bf 134} (2010), 189--195.


\bibitem{Po} {H. Poincar\'e}, Sur l'int\'egration alg\'ebrique
des \'equations diff\'erentielles du premier ordre et du premier
degr\'e, {\it Rendiconti del Circolo Matematico di Palermo} {\bf
5} (1891), 161--191.


\bibitem{Se} {A. Seidenberg},  Reduction of singularities of
the differential equation $Ady=Bdx$,  {\it Amer. J. Math.} {\bf
90} (1968), 248--269.

\bibitem{Sm} S. Smale, Mathematical problems for the next century, {\it Math. Intelligencer} {\bf 20} (1998), no. 2, 7--15.

\bibitem{Vi} {O. Viro},  From the sixteenth Hilbert problem to
tropical geometry, {\it Japan J. Math.} {\bf 3} (2008), 185--214.

\bibitem{Wi} {G. Wilson},  Hilbert's sixteenth problem, {\it Topology} {\bf 17}
(1978), 53--74.

\bibitem{Zh} {Xiang Zhang}, Invariant algebraic curves and rational
first integrals of holomorphic foliations in $\mathbb CP(2)$, {\it
Science in China Series A}  {\bf 46} (2003),  271--279.

\end{thebibliography}
\end{document}